\documentclass[sn-mathphys]{sn-jnl}
\usepackage{amssymb}
\usepackage{mathrsfs,amsmath,amssymb,bm,color,dsfont}
\usepackage{harpoon}
\usepackage{lineno,hyperref}
\usepackage{booktabs}
\usepackage{multirow}
\usepackage{subfigure}
\usepackage{ifpdf}
\usepackage{cases}
\usepackage{color}
\usepackage{threeparttable}
\usepackage{appendix}

\jyear{2021}%

\theoremstyle{thmstyleone}%
\newtheorem{theorem}{Theorem}
\theoremstyle{thmstyletwo}%
\newtheorem{definition}{Definition}
\newtheorem{lemma}{Lemma}

\newtheorem{remark}{Remark}

\theoremstyle{thmstylethree}%

\begin{document}

\title[Generalized Young Measure Solutions for a Class of Parabolic Equations]{Generalized Young Measure Solutions for a Class of Quasilinear Parabolic Equations with Linear Growth}


\author*[1]{\fnm{Jingfeng} \sur{Shao}}\email{sjfmath@hit.edu.cn}

\author[1]{\fnm{Zhichang} \sur{Guo}}\email{mathgzc@hit.edu.cn}
%
\author[1]{\fnm{Chao} \sur{Zhang}}\email{czhangmath@hit.edu.cn}

\affil[1]{\orgdiv{School of Mathematics}, \orgname{Harbin Institute of Technology}, \orgaddress{\street{West Da-Zhi Street}, \city{Harbin}, \postcode{150001}, \country{China}}}

%


\abstract{Using the generalized Young measure theory, we extend the theory of Young measure solutions to a class of quasilinear parabolic equations with linear growth, and introduce the concept of generalized Young measure solutions. We prove the existence and uniqueness of the generalized Young measure solutions. In addition, for the gradient flow of convex parabolic variational integral, we show that the generalized Young measure solutions are equivalent to the strong solutions.}




\maketitle

\section{Introduction}
Let $\Omega\subset \mathbb{R}^N$ be a bounded domain with
Lipschitz boundary $\partial \Omega$. Define $\Omega_T :=\Omega\times(0,T)$ and $\Gamma:=\partial\Omega \times (0,T)$. We study the following nonlinear
parabolic evolution problem of forward-backward type:
\begin{align}\tag{$\mathcal{P}$}\label{P}
\left\{\begin{array}{ll}
\dfrac{{\partial u}}{{\partial t}} = {\rm{div}}{\;\overrightarrow{q} \left( \nabla u \right)}\quad &\text{in }\Omega_T,
 \\u(x,t) = 0\quad &\text{on }\Gamma,
 \\u(x,0) = u_0(x)\quad &\text{in }\Omega.
\end{array}\right.
\end{align}
The initial data function $u_0$ is given in $BV(\Omega)\cap L^{\infty}(\Omega)$ with a trace of zero on $\partial \Omega$.
Here $\overrightarrow{q}:\mathbb{R}^N\rightarrow \mathbb{R}^N $ is a nonlinear, continuous, potential gradient function, satisfying that $\overrightarrow{q}=\nabla \Phi$, where $\Phi \in C^1 (\mathbb{R}^N)$ and $\Phi(0)=\inf\limits_{A\in \mathbb{R}^N} \Phi(A)$. Moreover, $\Phi$ and $\overrightarrow{q}$ satisfy the following structure conditions:
\begin{align}
{\left( {\lambda \vert A \vert - 1} \right)_ + } \leqslant \Phi (A) \leqslant \Lambda \vert A \vert + 1,\;\;\forall A\in\mathbb{R}^N
\end{align}
and
\begin{align}
\vert \overrightarrow{q}(A)\vert \leqslant \Lambda,\;\;\forall A\in\mathbb{R}^N,
\end{align}
for some $0<\lambda\leq \Lambda$, where $s_+=\max\{s,0\}$.
We would like to point out that no monotonicity requirements such as $(\overrightarrow{q}(A_1)-\overrightarrow{q}(A_2))\cdot(A_1-A_2)\geq 0$ for $A_1,A_2 \in \mathbb{R}^N$ are imposed.

Before stating our main results, let us mention some related results. Slemrod \cite{1991Dynamics} investigated the asymptotic behavior of measure valued solutions to the initial value problem for the nonlinear heat conduction equation
\begin{align*}
\dfrac{{\partial u}}{{\partial t}} = {\rm{div}}{\;\overrightarrow{q} \left( \nabla u \right)}\quad \text{in} \;\Omega_T,
\end{align*}
where
\[   {\lambda (\vert A \vert ^2 -1)}   \leqslant \Phi (A) \quad\text{and}\quad  \vert \overrightarrow{q}(A)\vert \leqslant \Lambda (1+ \vert A \vert^{\gamma}), \; 1\leq \gamma<2.    \]
Demoulini \cite{1996Young} further discussed the Young measure solutions for the case that
\begin{align*}
{\left( {\lambda \vert A \vert ^2 - 1} \right)_ + } \leqslant \Phi (A) \leqslant \Lambda \vert A \vert^2 + 1
\end{align*}
and
\begin{align*}
\vert \overrightarrow{q}(A)\vert \leqslant \Lambda \vert A \vert.
\end{align*}
With the similar framework as in \cite{1996Young}, Yin and Wang \cite{2003Young} extended the theory of Young measure solutions to the more general case like that
\begin{align*}
{\left( {\lambda \vert A \vert ^{1+\delta} - 1} \right)_ + } \leqslant \Phi (A) \leqslant \Lambda \vert A \vert ^{1+\delta} + 1
\end{align*}
and
\begin{align*}
\vert \overrightarrow{q}(A)\vert \leqslant \Lambda \vert A \vert^{\delta},
\end{align*}
where $0\leq \delta \leq 1$. Especially, for the limit case that $\delta=0$, by assuming that there exists a sequence $\{\Phi_\delta\}_{0<\delta\leq 1}$ such that $\{\overrightarrow{q}_\delta =\nabla \Phi_\delta\}_{0<\delta\leq 1}$ locally and uniformly converges to $\overrightarrow{q}$ in $\mathbb{R}^N$, and adding some conditions for $\{\Phi_\delta\}_{0<\delta\leq 1}$, they obtained the existence of a kind of biting Young measure solutions. However, the uniqueness is lost due to the failure of the equality of independence \eqref{idp} below.

When $\Phi$ is convex and $\delta=0$, there have been a lot of research activities in mathematical theory and engineering, see \cite{MR571133,MR1301176,2001Minimizing,MR3341130}. One possible approach to consider solutions is to use parabolic variational inequality. Lichnewsky and Temam \cite{1978Pseudo} considered the pseudo-solutions of the time-dependent minimal surface equation
$$\frac{{\partial u}}{{\partial t}} = {\text{div}}\;\left( {\frac{{\nabla u}}{{\sqrt {1{\text{ + }}\vert {\nabla u} \vert ^2} }}} \right). $$
More similar results can be found in \cite{MR1148668,MR2001660}.

Other possible approaches to define the generalized solutions is to apply the Anzellotti pairing \cite{Gabriele1983Pairings}. Andreu et al. \cite{2004Parabolic} investigated the following nonlinear problem
$$\frac{{\partial u}}{{\partial t}} = {\text{div}}\;\left( {\overrightarrow{a}(x,\nabla u)} \right),$$
where $\overrightarrow{a}(x,\xi)=\nabla_\xi f(x,\xi)$ and $f(x,\cdot)$ is a convex function with linear growth. By virtue of the Anzellotti pairing and nonlinear semigroup theory, they obtained the result of the existence and uniqueness of generalized solutions.

The aim of this paper is to study the existence and uniqueness of solutions for problem $(\mathcal{P})$. The novelties of this paper are as follows. First, inspired by the known theory \cite{2010Characterization,MR2885572}, we introduce a new framework of generalized Young measure solutions which can ensure the uniqueness of solutions. Second, our proof relies on some viscosity approximation instead of the usual $p$--Laplace type approximation as $p\to 1$. Third, based on the new definition of generalized Young measure solutions, we remark that we can not only deal with the nonlinear parabolic problem with linear growth, but also show the equivalence between the generalized Young measure solutions and the strong solutions (in Definition \ref{sms}).

This paper is organized as follows. In Section 2, we first present the mathematical preliminaries, some definitions and auxiliary lemmas. In Section 3, we prove the existence of the generalized Young measure solutions. In Section 4, we establish the uniqueness of the generalized Young measure solutions. Finally in Section 5, the equivalence between generalized Young measure solutions and strong solutions is investigated.

\section{Mathematical Preliminaries}

\subsection{The Auxiliary Problem and Young Measure Solutions}
By $\mathcal{M}(X)$ we denote the set of finite Radon measures on a Borel set $X$ and by the $M^+_1(X)$ its subset of probality measures. As usual $\mathcal{L}^N$ and $\mathcal{H}^{N-1}$ represent the $N$-dimensional Lebesgue measure and $(N-1)$-dimensional Hausdorff measure in $\mathbb{R}^N$, respectively.
We denote $\mathcal{M}_+(X)$ the positive Radon measures on $X$. For all $\lambda_1,\lambda_2\in \mathcal{M}(X)$,
\[{\lambda _1} \geqslant {\lambda _2}\;\text{ in }\;\mathcal{M}(X) \Leftrightarrow {\lambda _1} - {\lambda _2} \in \mathcal{M}_+(X).\]

Denote ${C_0}(\mathbb{R}^N)$ as the closure of continuous functions on $\mathbb{R}^N$ with compact support. The dual of ${C_0}(\mathbb{R}^N)$ can be identified with the Radon measures space $\mathcal{M}(\mathbb{R}^N)$ via the pairing
\[\left\langle {\nu ,f} \right\rangle  = \int\limits_{\mathbb{R}^N} f d{\nu },\quad   \forall \nu\in \mathcal{M}(\mathbb{R}^N),\;f \in {C_0}(\mathbb{R}^N)  .\]
Let the set $D \subset {\mathbb{R}^n}$ be a measurable set with finite measure. A map $\nu :D \to \mathcal{M} (\mathbb{R}^N)$ is called $weakly^\ast$ measurable if the functions $x \mapsto \int\limits_{\mathbb{R}^N} f d{\nu _x}$ are measurable for all $f \in {C_0}(\mathbb{R}^N)$, where ${\nu _x} = \nu (x)$.
Young measures on a bounded domain $\Omega\subseteq \mathbb{R}^N$ are $weakly^\ast$ measurable mappings $x\mapsto \nu_x,\;x\in\Omega$, with $\nu_x: \Omega\to \mathcal{M}^+_1(\mathbb{R}^N)$.

Next, the following lemma is a foundamental theorem on Young measures.
\begin{lemma}[\cite{MR1036070}]\label{lemma1} Let ${z^j}:D \to \mathbb{R}^N(j \ge 1)$ be a sequence of measurable functions. Then there exists a subsequence $\{ {z^{{j_k}}}\} _{k = 1}^\infty $ and a $weakly^\ast$ measurable map $\nu :D \to \mathcal{M} (\mathbb{R}^N)$ such that the following holds:
\begin{itemize}
  \item [\rm{(i)}] $\nu (x) \ge 0,\;{\left\| {\nu (x)} \right\|_{\mathcal{M}(\mathbb{R}^N)}} = \int\limits_{\mathbb{R}^N}  d{\nu _x} \le 1\quad\text{a.e.}\;x \in D$;\\
  \item [\rm{(ii)}] For all $f \in {C_0}(\mathbb{R}^N),\;f({z^{{j_k}}}(x))$ converges ${weakly^\ast}$ to $\int\limits_{\mathbb{R}^N} f d{\nu _x}$ in ${L^\infty }(D)$;\\
  \item [\rm{(iii)}] Furthermore, one has
             \[ {\left\| {\nu(x)} \right\|_{\mathcal{M}(\mathbb{R}^N)}} =1 \quad\text{a.e.} \;x\in D,\]
if and only if the sequence does not escape to infinity, namely,
\[\mathop {\lim }\limits_{M \to \infty }\mathop {\sup }\limits_{k \ge 1} \;\mathcal{L}^N(\{ x \in D:\vert {{z^{{j_k}}}} \vert \ge M\})  = 0.\]
\end{itemize}
\end{lemma}

\begin{definition} The map $\nu :D \to \mathcal{M}(\mathbb{R}^N)$ in Lemma \ref{lemma1} is called the Young measure on $\mathbb{R}^N$ generated by the sequence $\{ {z^{{j_k}}}\} _{k = 1}^\infty $.

For $p\geq 1$, define
\[\mathscr{E}_0^p(\mathbb{R}^N) = \left\{ {\varphi  \in C(\mathbb{R}^N):\mathop {\lim }\limits_{\vert A \vert \to \infty } \frac{{\vert {\varphi (A)} \vert}}{{1 + {{\vert A \vert}^p}}}\;\;exists} \right\}.\]
The space $\mathscr{E}_0^p(\mathbb{R}^N)$ is a separable Banach space with the norm
\[{\left\| \phi  \right\|_{\mathscr{E}_0^p}} = \mathop {\sup }\limits_{A \in \mathbb{R}^N} \frac{{\vert {\phi (A)} \vert}}{{1 + {{\vert A \vert}^p}}}.\]

Define
\[{\mathscr{E}^p}(\mathbb{R}^N) = \left\{ {\varphi  \in C(\mathbb{R}^N):\mathop {\sup }\limits_{A \in \mathbb{R}^N}  \frac{{\vert {\varphi (A)} \vert}}{{1 + {{\vert A \vert}^p}}}\; <  + \infty } \right\},\]
which is an inseparable space in the above norm.
\end{definition}

Let $\Phi ^{**}$ be the convexification of $\Phi$, namely,
\begin{align*}
{\Phi ^{**}}(A) = \sup \left\{ {f(A):f \leq \Phi ,\;f\;\text{is convex}} \right\}.
\end{align*}
Since $\Phi \in C^1 (\mathbb{R}^N)$, then $\Phi^{**}$ is convex and in $C^1 (\mathbb{R}^N)$ (see \cite{KIRCHHEIM2001725}).

\begin{lemma}[Theorem 2.35, \cite{2007Direct}] \label{envelope}
Let $f:\mathbb{R}^N \to \mathbb{R}\cup \{+\infty\} $ and, for every $x\in \mathbb{R}^N$,
\begin{align*}
Cf(x)=\sup\left\{g(x):g\leq f {\text{ and }} g {\text{ is convex}}\right\}.
\end{align*}
Assume that $Cf>-\infty$. Then
\begin{align*}
Cf(x)={\rm inf} \left\{ \sum\limits_{i = 1}^{N + 1} {{\alpha _i}f({x_i})} :\sum\limits_{i = 1}^{N + 1} {{\alpha _i}{x_i}}  = x,\;{\alpha _i} \geqslant 0\;\;{\text{with}}\;\sum\limits_{i = 1}^{N + 1} {{\alpha _i}}  = 1\right\}.
\end{align*}
\end{lemma}

In the rest of this paper, $I$ represents the identity operator in $\mathbb{R}^N$.

Set $\overrightarrow{p}=\nabla \Phi ^{**}$.
We note that $\overrightarrow{q}=\overrightarrow{p}$ on the set $\left\{A\in \mathbb{R}^N:\Phi(A)=\Phi^{**}(A)\right\}$. Thus, $\Phi^{**}$ and $\overrightarrow{p}$ satisfy the same structure conditions as $\Phi$ and $\overrightarrow{q}$, respectively.

We consider the following auxiliary problem
\begin{equation}\tag{$\mathcal{P_\varepsilon}$}\label{Pep}
\left\{\begin{array}{ll}
\dfrac{{\partial u}}{{\partial t}} ={\rm{div}}({\;\overrightarrow{q}_\varepsilon \left( \nabla u \right)}):= {\rm{div}}({\;\overrightarrow{q} \left( \nabla u \right)+\varepsilon\nabla u})\quad &\text{in }\Omega_T,
 \\u(x,t) = 0\quad &\text{on }\Gamma,
 \\u(x,0) = u^{\varepsilon}_0(x)\quad &\text{in }\Omega,
\end{array}\right.
\end{equation}
where $ u^{\varepsilon}_0(x)\in H^1_0(\Omega)$ and it satisfies that when $\varepsilon \rightarrow 0$
\begin{numcases}{}
\left\| {u_0^\varepsilon (x)} \right\|_{BV(\Omega )} \to {\left\| {u_0(x)} \right\|_{BV(\Omega )}},
\\u_0^\varepsilon (x) \to {u_0}(x)\quad\text{in}\;{L^2}(\Omega ),
\\{\left\| {u_0^\varepsilon (x)} \right\|_{{L^\infty }(\Omega )}} \leqslant {\left\| {{u_0}(x)} \right\|_{{L^\infty }(\Omega )}},
\\\sqrt {{\varepsilon _n}} {\left\| {u_0^{{\varepsilon _n}}(x)} \right\|_{{H^1}(\Omega )}} \to 0,\;\text{for a subsequence of}\; \{\varepsilon\}.
\end{numcases}

\begin{definition}  A Young measure valued solution of $(\mathcal{P_\varepsilon})$ is a pair $(u^\varepsilon,\nu ^\varepsilon)$, where
$u^\varepsilon \in {L^{\infty}}(0,T;H_0^1(\Omega ))$,
$\dfrac{{\partial u^\varepsilon}}{{\partial t}} \in L^2(\Omega_T)$, and $\nu^\varepsilon  = {\left( {{\nu ^\varepsilon _{x,t}}} \right)_{\left( {x,t} \right) \in {\Omega_T}}}$ is a parametrized family of probability measures on $\mathbb{R}^N$ such that
\begin{align}
\iint_{{\Omega_T}} { {\left( {\left\langle {\nu^\varepsilon ,\overrightarrow{q}_\varepsilon} \right\rangle  \cdot \nabla \zeta  + \frac{{\partial u ^\varepsilon}}{{\partial t}}\zeta } \right)} } dxdt = 0,\;\forall \zeta  \in H^1_0({\Omega_T})&,\label{au1}\\
\nabla u^\varepsilon(x,t) = \left\langle {{\nu _{x,t}},I} \right\rangle\quad \text{a.e.}\;(x,t) \in {\Omega_T} \label{au2}&,
\end{align}
and
\begin{align}
u^\varepsilon(x,t) &= 0,\;(x,t) \in \partial \Omega  \times [0,T] \label{au3},
\\u^\varepsilon(x,0) &= {u^\varepsilon_0}(x),\;x \in \Omega \label{au4},
\end{align}
in the sense of trace. Furthermore,
\begin{align}
{\mathrm{supp}}\;\nu^\varepsilon  \subseteq \left\{ {\Phi  + \frac{\varepsilon}{2}{{\vert I \vert}^2} = {{\left( {\Phi  + \frac{\varepsilon}{2}{{\vert I \vert}^2}} \right)}^{**}}} \right\},
\end{align}
and the equality of independence
\begin{align}
\left\langle {{\nu^\varepsilon _{x,t}},\overrightarrow{q}_\varepsilon  \cdot I} \right\rangle  = \left\langle {{\nu^\varepsilon _{x,t}},\overrightarrow q_\varepsilon } \right\rangle  \cdot \left\langle {{\nu^\varepsilon _{x,t}},I} \right\rangle \quad\text{a.e.}\;(x,t)\in{\Omega _T} \label{idp}
\end{align}
holds, i.e., $\overrightarrow{q}_\varepsilon$  and $\nabla u^\varepsilon$ are independent with respect to the Young measure $\nu^\varepsilon$.
\end{definition}

According to Remark 3.1 in \cite{2003Young}, the following lemma holds.
\begin{lemma} \label{lemma2}
The problem $(\mathcal{P}_\varepsilon)$ admits a Young measure solution $u^\varepsilon$. And there exists a constant $M$ depending only on ${\left\| {{u_0}(x)} \right\|_{{BV}(\Omega )}},{\left\| {{u_0}(x)} \right\|_{{L^\infty }(\Omega )}},\lambda ,\Lambda $, and $\mathcal{L}^N (\Omega)$ but independent of $\varepsilon$ and $T$ such that
\begin{numcases}{}
{\left\| {{u^\varepsilon }} \right\|_{{L^\infty }({\Omega _T})}} \leqslant M,
\\{\left\| {{u^\varepsilon }} \right\|_{{L^\infty }(0,T;{W^{1,1}_0}(\Omega ))}} \leqslant M,
\\{\left\| {{u^\varepsilon }} \right\|_{{L^\infty }(0,T;{L^2}(\Omega ))}} \leqslant M,
\\{\left\| {\frac{{\partial {u^\varepsilon }}}{{\partial t}}} \right\|_{{L^2}\left( {{\Omega _T}} \right)}} \leqslant M,
\\\sqrt {{\varepsilon _n}} {\left\| {{u^{{\varepsilon _n}}}} \right\|_{L^\infty(0,T;{H^1}(\Omega ))}} \leqslant M.
\end{numcases}
\end{lemma}

\subsection{Pairings between Measures and Bounded Functions}

A function $u\in L^1(\Omega)$ is called a function of bounded variation if its distributional derivative $Du$ is a $\mathbb{R}^N$-valued Radon measure with finite total variation in $\Omega$. The vector space of functions of bounded variation in $\Omega$ is denoted by $BV(\Omega)$. The space $BV(\Omega)$ is a non-reflexive Banach space under the norm $\| u\|_{BV}:=\| u\|_{L^1} + \vert Du \vert (\Omega)$, where $\vert Du \vert$ is the total variation measure of $Du$.

For $u\in BV(\Omega)$, the gradient $Du$ is a Radon measure that decomposes into its absolutely continuous and singular parts
\[ Du=D^a u +D^s u .\]
Then $D^a u = \nabla u \mathcal{L}^N$, where $\nabla u$ is the Radon-Nikod\'{y}m derivative of the measure $Du$ with respect to the Lebesgue measure $\mathcal{L}^N$. There is also the polar decomposition $D^s u = \frac{{{D^s}u}}{{\vert {D^s}u \vert}}\vert {{D^s}u} \vert $, where $\frac{{{D^s}u}}{{\vert {D^s}u \vert}}  \in {L^1}\left( {\Omega ,\vert {{D^s}u} \vert ;\partial {\mathbb{B}^N}} \right)$ is the Radon-Nikod\'{y}m derivative of $D^s u$ with respect to its total variation measure $\vert D^s u \vert $.

Set $\tilde{\Omega}$ as a bounded smooth domain satisfying $\Omega \subset\subset \tilde{\Omega}$. For $v \in BV(\Omega)$, we denote $\tilde{v}$ by
\[\tilde v(x): =
\left\{ \begin{gathered}
  v(x),\;x \in \Omega ,
  \\0, \;x \in \tilde{\Omega}\backslash \overline{\Omega}.
\end{gathered}  \right.\]

\begin{lemma}[Theorem B.3, \cite{2004Parabolic}]
Assume that $u\in BV(\Omega)$. There exists a sequence of functions $u_i \in W^{1,1}(\Omega)\cap C^\infty (\Omega)$ such that
\begin{itemize}
\item [\rm{(i)}] $u_i \to u$ in $L^1(\Omega)$;
\item [\rm{(ii)}] $\vert Du_i\vert(\Omega) \to \vert Du \vert(\Omega)$;
\item [\rm{(iii)}] $u_i \vert_{\partial \Omega} =u \vert_{\partial \Omega}$ for all $i$;
\end{itemize}
Moreover,
\begin{itemize}
\item [\rm{(iv)}] if $u \in BV(\Omega)\cap L^q (\Omega)$, $q<\infty$, we can find functions $u_i$ such that $u_i \in L^q(\Omega)$ and $u_i \to u$ in $L^q(\Omega)$;
\item [\rm{(v)}] if $u \in BV(\Omega)\cap L^\infty (\Omega)$, we can find $u_i$ such that $\|u_i\|_\infty \leq \|u\|_\infty$ and $u_i \to u$ weakly* in $L^\infty(\Omega)$.
\end{itemize}
\end{lemma}

It is well known that summability conditions on the divergence of a vector field $z$ in $\Omega$ yield trace properties for the normal component of $z$ on $\partial \Omega$. As in \cite{Gabriele1983Pairings}, we define a function $[z,\vec{n}]\in L^\infty (\partial\Omega)$ which is associated to any vector field $z\in L^\infty (\Omega, \mathbb{R}^N)$ such that
${\rm{div}} (z)$ is a bounded measures in $\Omega$.

Assume that $\Omega$ is an open bounded set in $\mathbb{R}^N$ with $\partial \Omega$ Lipschitz, $N\geq 2$, and $1\leq p \leq N$, $\frac {N}{N-1}\leq q \leq \infty$. Since $\partial \Omega$ is Lipschitz, the outer unit normal $\vec{n}$ exists $\mathcal{H}^{N-1}-\mathrm{a.e.}$ on $\partial \Omega$. We shall consider the following spaces:
\begin{align*}
&BV(\Omega)_q:=BV(\Omega)\cap L^q(\Omega),
\\&BV(\Omega)_c:=BV(\Omega)\cap L^\infty(\Omega)\cap C(\Omega),
\\&X(\Omega)_p:=\{z\in L^\infty(\Omega,\mathbb{R}^N):{\rm{div}} (z)\in L^p(\Omega)\},
\\&X(\Omega)_\mu :=\{z\in L^\infty(\Omega,\mathbb{R}^N):{\rm{div}} (z)\in \mathcal{M}(\Omega)\}.
\end{align*}

\begin{lemma}[Theorem C.2, \cite{2004Parabolic}]\label{X1}
Assume that $\Omega \subseteq \mathbb{R}^N$ is an open bounded set with Lipschitz boundary $\partial \Omega$. Then there exists a bilinear map $\langle z,u\rangle_{\partial \Omega} : X(\Omega)_\mu \times BV(\Omega)_c \to \mathbb{R}$ such that
\begin{align}
\langle a, u\rangle_{\partial \Omega}=\int\limits_{\partial \Omega } {u(x)z(x) \cdot \vec n\;d{\mathcal{H}^{N - 1}}}  \;\;\text{if}\; z\in C^1(\overline{\Omega},\mathbb{R}^N),
\end{align}
\begin{align}
\vert\langle a, u\rangle_{\partial \Omega}\vert\leq \| z\|_\infty \int\limits_{\partial \Omega } \vert u(x)\vert d {\mathcal{H}^{N - 1}}.
\end{align}
\end{lemma}

\begin{lemma}[Theorem C.3, \cite{2004Parabolic}]
Assume that $\Omega \subseteq \mathbb{R}^N$ is an open bounded set with Lipschitz boundary $\partial \Omega$. Then there exists a linear operator $\gamma: X(\Omega)_\mu \to L^\infty(\partial \Omega)$ such that
\begin{align}
\| \gamma\|_\infty &\leq \| z\|_\infty,
\\ \langle a, u\rangle_{\partial \Omega}&=\int\limits_{\partial \Omega } {u(x)\gamma (z)(x) \cdot \vec n\;d{\mathcal{H}^{N - 1}}}  \;\;\text{ for all}\; u\in BV(\Omega)_c,
\\ \gamma (z)(x)&=z(x)\cdot \vec{n} \;\;\; \text{ for all}\; x\in \partial \Omega \;\;\text{if}\;\; z\in C^1(\overline{\Omega},\mathbb{R}^N).
\end{align}
The function $\gamma(z)$ is a weakly defined trace on $\partial \Omega$ of the normal component of $z$. We denote $\gamma (z)$ by $[z,\vec{n}]$.
\end{lemma}

In the sequel we shall consider pairs $(z,u)$ such that one of the following conditions holds
\begin{numcases}{}
a)\; z \in X(\Omega)_p,\;u\in BV(\Omega)_{p'} \;\text{and}\; 1<p \leq N;\nonumber
\\b)\; z \in X(\Omega)_1,\;u\in BV(\Omega)_{\infty}; \label{pair}
\\c) \; z \in X(\Omega)_\mu,\;u\in BV(\Omega)_c.\nonumber
\end{numcases}

\begin{definition}
Let $(z,u)$ satisfy one of the conditions \eqref{pair}. Then we define a functional $(z,Du):\mathcal{D}(\Omega) \to \mathbb{R}$ as
\begin{align*}
\langle (z,Du),\phi \rangle := - \int\limits_{\Omega } {u(x)\phi\; {\rm{div}}(z)dx}- \int\limits_{\Omega } {u(x)z\cdot \nabla \phi dx}, \;\;\forall\; \phi \in \mathcal{D}(\Omega),
\end{align*}
where $\mathcal{D}(\Omega)$ represents the space of all infinitely differentiable functions with compact support on $\Omega$.
\end{definition}

\begin{lemma}[Theorem C.6, \cite{2004Parabolic}]
For any functions $\phi \in \mathcal{D}(V)$ and for any open sets $V\subseteq \Omega$, it follows that
\begin{align}
\vert \langle (z,Du),\phi \rangle \vert \leq {\sup}\; \|\phi \|_\infty \| z\|_{L^\infty(V)} \vert Du \vert(V).
\end{align}
Thus, $(z,Du)$ is a Radon measure in $\Omega$.
\end{lemma}

Next, we give the Green's formula related to the function $[z,\vec{n}]$ and the measure $(z,Du)$.

\begin{lemma}[Theorem C.9, \cite{2004Parabolic}]\label{green}
Assume that $\Omega \subseteq \mathbb{R}^N$ is an open bounded set with Lipschitz boundary $\partial \Omega$. If $z$ and $u$ ensure that one of the conditions \eqref{pair} holds, then we have
\begin{align}
\int\limits_{\Omega } {u\; {\rm{div}}(z)dx} + \int\limits_{\Omega } {(z,Du)} = \int\limits_{\partial\Omega } {[z,\vec{n}]\;u \;d \mathcal{H}^{N-1}}.
\end{align}
\end{lemma}
More details about pairings $\langle z, u\rangle_{\partial\Omega}$ and $(z,Du)$ can be found in \cite{Gabriele1983Pairings}.


\subsection{The Generalized Young Measure Solutions}

This subsection gives a brief overview of the basic theory of generalized Young measures, which will be used later. We refer to \cite{2010Characterization,MR2885572} for further information about the generalized Young measures.

First, we need a suitable class of integrands. Let $\mathcal{E}(\Omega; \mathbb{R}^N)$ be the set of all $f\in C(\overline{\Omega}\times \mathbb{R}^N)$ such that
\begin{align}
(Tf)(x,\hat A): = (1 - \vert\hat A \vert)f\left( {x,\frac{{\hat A}}{{1 - \vert \hat A\vert}}} \right),\;x \in \overline \Omega,\;\hat A \in {\mathbb{B}^N},
\end{align}
extends into a continuous function $Tf \in C(\overline {\Omega  \times {\mathbb{B}^N}})$, where $\mathbb{B}^N$ denotes the open unit ball in $\mathbb{R}^N$. In particular, it implies that $f$ has linear growth at infinity, i.e., there exists a constant $M> 0$ (in fact, $M  ={\left\| f \right\|_{\mathcal{E}(\Omega ;{\mathbb{R}^N})}}: = {\left\| {Tf} \right\|_{\infty ,\overline {\Omega  \times {\mathbb{B}^N}} }}$ ),
\[\vert {f(x,A)} \vert \leqslant M(1 + \vert A \vert)\;\;\text{for all}\;(x,A) \in \overline \Omega   \times {\mathbb{R}^N}.\]
For all $f\in \mathcal{E}(\Omega; \mathbb{R}^N)$ the recession function
\begin{align}
{f^\infty }(x,A): = \mathop {\lim }\limits_{\substack{
 A' \to A \\ x' \to x\\t \to \infty}}
\frac{{f(x',tA')}}{t},\;x \in \overline \Omega,\;A \in {\mathbb{R}^N},
\end{align}
exists as a continuous function. Sometimes this notion of a recession function is too strong and so for any function $g\in C(\mathbb{R}^N)$ with linear growth at infinity, we define the generalized recession function
\begin{align}
{g^\sharp}(A): = \mathop {\lim }\limits_{\substack{
 A' \to A \\t \to \infty}}
\frac{{g(tA')}}{t},\;A \in {\mathbb{R}^N}.
\end{align}
We remark that for both flavors of recession function one can drop the additional sequence $A'\to A$ if the functional is Lipschitz continuous (see \cite{2014A}).
\begin{definition} A generalized Young measure with target space $\mathbb{R}^N$ is a triple $\lambda=(\nu_x,\lambda_\nu ,{\nu}^\infty _x)$ comprising
\begin{itemize}
\item [\rm{(i)}] a parametrized family of probability measures $(\nu_x)_{x\in\Omega} \subseteq \mathcal{M}^+_1(\mathbb{R}^N)$;
\item [\rm{(ii)}] a positive finite measure $\lambda_\nu \in \mathcal{M}_+(\overline{\Omega})$;
\item [\rm{(iii)}] a parametrized family of probability measures $(\nu^\infty _x)_{x\in \overline{\Omega}}\subseteq \mathcal{M}^+_1(\partial\mathbb{B}^N)$ ($\partial\mathbb{B}^N$ denotes the unit sphere in $\mathbb{R}^N$);
\item [\rm{(iv)}] the map $x\mapsto \nu_x$ is weakly* measurable with respect to $\mathcal{L}^N$, i.e., the function $x\mapsto \langle \nu_x , f(x,\cdot)\rangle$ is $\mathcal{L}^N$-measurable for all bounded Borel functions $f:\Omega\times \mathbb{R}^N \to \mathbb{R}$;

\item [\rm{(v)}] The map $x\mapsto \nu^\infty _x$ is weakly* measurable with respect to $\lambda_\nu$;

\item [\rm{(vi)}] $x\mapsto \langle \nu_x, \vert \cdot \vert  \rangle \in L^1 (\Omega)$.

\end{itemize}

\end{definition}

By $Y(\Omega;\mathbb{R}^N)$ we denote the set of all such generalized Young measures. The parametrized measure $(\nu_x)$ is called the \textbf{oscillation measure}, the measure $\lambda_\nu$
is the \textbf{concentration measure}, and $(\nu^\infty _x)$ is the \textbf{concentration-angle measure}.

The duality  pairing $\langle\langle \lambda,f \rangle\rangle$ for $\lambda \in Y(\Omega;\mathbb{R}^N)$ and $f\in \mathcal{E}(\Omega; \mathbb{R}^N)$ is defined via
\begin{align*}
  \left\langle {\left\langle {\lambda ,f} \right\rangle } \right\rangle &: = \int\limits_\Omega  {\left\langle {{\nu _x},f(x, \cdot )} \right\rangle dx}  + \int\limits_{\overline \Omega  } {\left\langle {\nu _x^\infty ,{f^\infty }(x, \cdot )} \right\rangle d{\lambda _\nu }(x)}
 \\&  = \int\limits_\Omega  {\int\limits_{{\mathbb{R}^N}} {f(x,A)d{\nu _x}\left( A \right)} dx}  + \int\limits_{\overline \Omega  } {\int\limits_{\partial {\mathbb{B}^N}} {{f^\infty }(x, \cdot )d\nu _x^\infty \left( A \right)} d{\lambda _\nu }(x)}.
\end{align*}

The space $Y(\Omega; \mathbb{R}^N)$ of Young measures can be considered as a part of the dual space $\mathcal{E}(\Omega;\mathbb{R}^N)^*$ (the inclusion is strict since, for instance, $f\mapsto r\int\limits_\Omega  {f(x, A_0 )dx} $ lies in $\mathcal{E}(\Omega;\mathbb{R}^N)^*\setminus Y(\Omega; \mathbb{R}^N)$ whenever $r\ne 1$ and $A_0\in \mathbb{R}^N$). This embedding gives rise to a weak* topology on $Y(\Omega; \mathbb{R}^N)$ and so we say that $(\lambda_j)\subseteq Y(\Omega; \mathbb{R}^N)$ (where $\lambda_j:=(\nu^j _x,\lambda_{\nu_j},{\nu}_x^{\infty,j}$) weakly* converges to $\lambda \in Y(\Omega; \mathbb{R}^N)$, in symbols
$\lambda_j \stackrel{*}{\rightharpoondown} \lambda$, if $\langle\langle \lambda_j, f\rangle\rangle \to \langle\langle \lambda, f\rangle\rangle$ for all $f\in\mathcal{E}(\Omega;\mathbb{R}^N)$. The set $Y(\Omega; \mathbb{R}^N)$ is topologically weakly*-closed in $\mathcal{E}(\Omega;\mathbb{R}^N)^*$.

The main compactness result in the space $Y(\Omega; \mathbb{R}^N)$ is listed as follows.
\begin{lemma}[Corollary 2, \cite{2010Characterization}]
Let $(\lambda_j)\subseteq Y(\Omega; \mathbb{R}^N)$ be a sequence such that
\begin{itemize}
\item [\rm{(i)}] the functions $x\mapsto \langle \nu_x^j, \vert\cdot\vert\rangle$ are uniformly bounded in $L^1(\Omega)$;
\item [\rm{(ii)}] the sequence $(\lambda_{\nu_j}(\overline{\Omega}))$ is uniformly bounded.
\end{itemize}
Then, $(\lambda_j)$ is weakly* sequentially relatively compact in $Y(\Omega; \mathbb{R}^N)$, i.e., there exists $\lambda \in Y(\Omega; \mathbb{R}^N)$ and a subsequence of $(\lambda_j)$ (not relabeled) such that $\lambda_j \stackrel{*}{\rightharpoondown} \lambda $.
\end{lemma}

Next, we define the set $GY(\Omega; \mathbb{R}^N)$ of generalized gradient Young measures as the collection of the generalized Young measures $\lambda \in Y(\Omega; \mathbb{R}^N)$ with the property that there exists a norm-bounded sequence $(u_j)\subseteq BV(\Omega)$ such that the sequence $(Du_j)$ generates $\lambda$, denoted as $Du_j \stackrel{Y}{\to} \lambda$, meaning that
\[\int\limits_\Omega  {f(x,\nabla {u_j}(x))dx}  + \int\limits_{\overline \Omega  } {{f^\infty }(x,\frac{{{D^s}{u_j}}}{{\left\vert {{D^s}{u_j}} \right\vert}})d} \left\vert {{D^s}{u_j}} \right\vert \to \left\langle {\left\langle {\lambda ,f} \right\rangle } \right\rangle \]
for all $f\in \mathcal{E}(\Omega;\mathbb{R}^N)$.

The following lemma can been found in \cite{2014A}, see also \cite{2010Characterization}.
\begin{lemma}[Theorem 9, \cite{2010Characterization}] \label{charact}
Let $\lambda \in Y(\Omega;\mathbb{R}^N)$ be a generalized Young measure with $\lambda_\nu (\partial \Omega)=0$.
Then $\lambda \in GY(\Omega; \mathbb{R}^N)$, if and only if there exists $u\in BV(\Omega)$ with
\begin{align*}
Du =\left\langle {\nu,I} \right\rangle {\mathcal{L}^N}\llcorner\Omega + \left\langle {{\nu^\infty },I} \right\rangle {\lambda _\nu }\llcorner \Omega  \quad\text{ in }\;\mathcal{M}(\Omega),
\end{align*}
and for all quasiconvex $g\in C(\mathbb{R}^N)$ with linear growth at infinity the following Jensen-type inequalities hold:
\begin{itemize}
\item [\rm{(i)}] \[g\left( {\left\langle {\nu,I} \right\rangle  + \left\langle {{\nu^\infty },I} \right\rangle \dfrac{{d{\lambda _\nu }}}{{d{\mathcal{L}^N}}}} \right) \leqslant \left\langle {\nu,g} \right\rangle  + \left\langle {{\nu^\infty },{g^\sharp}} \right\rangle \dfrac{{d{\lambda _\nu }}}{{d{\mathcal{L}^N}}},\]
    for $\mathcal{L}^N$-a.e. $x\in \Omega$, where $\dfrac{{d{\lambda _\nu }}}{{d{\mathcal{L}^N}}}$ is the Radon-Nikod\'{y}m derivative of the measure $\lambda$ with respect to the Lebesgue measure $\mathcal{L}^N$;
\item [\rm{(ii)}] \[{g^\sharp}\left( {\left\langle {{\nu^\infty },I} \right\rangle } \right) \leqslant \left\langle {{\nu^\infty },{g^\sharp}} \right\rangle \]
for $\lambda^s _\nu$-a.e. $x\in \Omega$, where $\lambda^s _\nu$ is the singular part of $\lambda$ with respect to the Lebesgue measure $\mathcal{L}^N$.
\end{itemize}
\end{lemma}

Given $u\in BV(\Omega)$, denote $\sigma_{Du} :=(\delta_{\nabla u}, \vert D^s u \vert, \delta _p)\in GY(\Omega;\mathbb{R}^N)$ and $p:=\frac{{{D^s}u}}{{\vert {{D^s}u} \vert }} \in {L^1}\left( {\Omega ,\vert {{D^s}u} \vert;\partial {\mathbb{B}^N}} \right)$. $\delta_A$ denotes the Dirac measure on $\mathbb{R}^N$ giving unit mass to the point $A\in \mathbb{R}^N$.

In the rest of this paper, $u^\Omega$ represents the trace of $u\in BV(\Omega)$ on $\partial \Omega$.

\begin{definition}  A generalized Young measure valued solution of $(\mathcal{P})$ is a pair $(u,\lambda)$, if
$u \in {L^{\infty}}(0,T;BV(\Omega )\cap L^2(\Omega))$, $\dfrac{\partial u}{{\partial t}} \in L^2(\Omega_T)$, $\lambda=(\nu_{x,t},\lambda_\nu,\nu^\infty_{x,t})_{x\in\Omega}\in Y(\Omega;\mathbb{R}^N)$, $\langle\langle \vert \cdot \vert, \lambda \rangle\rangle \in L^\infty (0,T)$,  and for almost all $t\in (0,T)$ they satisfy
\begin{align}
&u_t=\rm{div} \langle \nu,\overrightarrow{q} \rangle \quad\text{in}\; \mathcal{D}'(\Omega),
\\&Du = \left\langle {\nu,I} \right\rangle {\mathcal{L}^N}\llcorner\Omega + \left\langle {{\nu^\infty },I} \right\rangle {\lambda _\nu }\llcorner \Omega \quad\text{in}\;\mathcal{M}(\Omega),
\end{align}
and
\begin{align}
&\left\langle {{\nu^\infty },I} \right\rangle {\lambda _\nu }\llcorner \partial\Omega =\left( { - u^\Omega \otimes \vec{n}} \right){\mathcal{H}^{N-1}}\llcorner \partial\Omega  \quad\text{in}\;\mathcal{M}(\overline{\Omega}),
\\&u(x,0) = {u_0}(x),\;x \in \Omega.
\end{align}
Furthermore,
\begin{align}
{\mathrm{supp}}\;\nu  \subseteq \left\{\Phi = \Phi^{**} \right\},   \label{supp}
\end{align}
and the {\bf{(JF)}} inequality
\begin{align}
&\left\langle {\nu ,\overrightarrow q } \right\rangle  \cdot \left(  \left\langle {\nu,I} \right\rangle {\mathcal{L}^N}\llcorner\Omega + \left\langle {{\nu^\infty },I} \right\rangle {\lambda _\nu }\llcorner \overline{\Omega}\right) \nonumber
\\\geq & \left\langle {\nu,\overrightarrow{q}\cdot I} \right\rangle {\mathcal{L}^N}\llcorner\Omega + \left\langle {{\nu^\infty },(\overrightarrow{q}\cdot I)^\infty} \right\rangle {\lambda _\nu }\llcorner \overline{\Omega}   \label{jf}
\end{align}
holds in $\mathcal{M}(\overline{\Omega})$.
\end{definition}

\begin{remark}
If $(u,\lambda)$ is a generalized Young measure solution, and for almost all $t\in(0,T)$, $\lambda \in GY(\Omega;\mathbb{R}^N)$, we call $(u,\lambda)$ is a generalized gradient Young measure solution.
\end{remark}


\section{Existence Results}
In this section, we are ready to prove the existence of generalized Young measure solutions of $(\mathcal{P})$.

We first introduce the following assumptions {\bf{(SH)}}:
\\$\rm{(SH_1)}$ $\Phi \in C^1(\mathbb{R}^N)$, $\Phi(0)=\inf\limits_{A\in \mathbb{R}^N} \Phi(A)$ and $\Phi \in \mathcal{E}(\Omega ; \mathbb{R}^N)$.
\\$\rm(SH_2)$ $\Phi$ and $\overrightarrow{q}=\nabla \Phi$ satisfy the structure conditions:
\begin{align}
{\left( {\lambda \vert A \vert - 1} \right)_ + } \leqslant \Phi (A) \leqslant \Lambda \vert A \vert + 1,\;\forall A\in\mathbb{R}^N,
\end{align}
and
\begin{align}
\vert \overrightarrow{q}(A)\vert\leqslant \Lambda,\;\;\forall A\in\mathbb{R}^N,
\end{align}
for some $0<\lambda\leq \Lambda$.
\\$\rm (SH_3)$ $\overrightarrow q  \cdot I \in \mathcal{E}(\Omega ; \mathbb{R}^N)$ and $(\overrightarrow q  \cdot I)^\infty={\Phi}^\infty=(\Phi^{**})^\infty$.
\\$\rm (SH_4)$ $\forall \eta \in \partial \mathbb{B}^N$, $\forall \xi \in \left\{ {A \in {R^N}\mid {{\Phi ^{**}}(A) = \Phi (A)} } \right\}$,
\[\overrightarrow q (\xi ) \cdot \eta  \leqslant {\left( {\overrightarrow q  \cdot I} \right)^\infty }(\eta ).\]

\begin{theorem}\label{exist}
Under the assumptions of {\bf{(SH)}}. Given $u_0\in BV(\Omega)\cap L^{\infty}(\Omega)$ with a trace of zero on $\partial \Omega$, there exists a generalized Young measure solution of $(\mathcal{P})$ for every $T>0$.
\end{theorem}
{\bf{Proof}}: According to Lemma \ref{lemma2}, for almost all $t\in(0,T)$, we can extract from $\{u^{\varepsilon_n}\}$ and $\{\nu^{\varepsilon_n}\}$, a subsequence (still labeled by $\{u^{\varepsilon_n}\}$ and $\{\nu^{\varepsilon_n}\}$) such that
\begin{numcases}{}
{u^{{\varepsilon _n}}} \to u\quad\text{in}\;{L^2}(\Omega ), \label{cve1}
\\{u^{{\varepsilon _n}}} \stackrel{*}{\rightharpoonup} u\quad\text{in}\;{L^\infty}(\Omega ),
\\u_t^{{\varepsilon _n}} {\rightharpoonup} {u_t}\quad\text{in}\;{L^2}(\Omega ), \label{cve2}
\\ \nabla \tilde{u}^{\varepsilon_n}{\mathcal{L}^N}\llcorner\tilde{\Omega} \stackrel{*}{\rightharpoondown} D\tilde{u} \quad\text{in}\;\mathcal{M}(\tilde{\Omega}),
\\{\varepsilon _n}\vert {\nabla {u^{{\varepsilon _n}}}} \vert \to 0\quad\text{in}\;{L^2}(\Omega ), \label{cve3}
\end{numcases}
and
\begin{numcases}{}
\left\langle {{\nu ^{{\varepsilon _n}}},\overrightarrow q } \right\rangle \stackrel{*}{\rightharpoonup} \left\langle {\nu ,\overrightarrow q } \right\rangle \quad\text{in}\;{L^\infty }(\Omega;\mathbb{R}^N), \label{ve1}
\\\left\langle {{\nu ^{{\varepsilon _n}}},\overrightarrow q  \cdot I} \right\rangle {\mathcal{L}^N}\llcorner\Omega \stackrel{*}{\rightharpoondown} \left\langle {\nu,\overrightarrow q  \cdot I} \right\rangle {\mathcal{L}^N}\llcorner\Omega + \nonumber
\\\left\langle {{\nu^\infty },{{(\overrightarrow q  \cdot I)}^\infty }} \right\rangle {\lambda _\nu }\quad\text{in}\;\mathcal{M}(\overline{\Omega}), \label{ve2}
\\{\lambda ^{{\varepsilon _n}}} \stackrel{*}{\rightharpoondown}  \lambda \quad\text{in}\;Y\left( {\Omega ;{\mathbb{R}^N}} \right),\label{ve3}
\end{numcases}
where $\lambda ^{\varepsilon _n}:=(\nu^{\varepsilon_n},0,\delta_0)$ and  $\lambda=(\nu,\lambda_\nu,\nu^\infty)$.

We divide the proof into the following steps.

{\textbf{Step 1}}. For all $\varphi \in \mathcal{D}(\Omega)$,
\[\int\limits_\Omega  {{\varepsilon _n}\left\langle {{\nu ^{{\varepsilon _n}}},I} \right\rangle  \cdot \nabla \varphi dx}  = \int\limits_\Omega  {\nabla {u^{{\varepsilon _n}}} \cdot \nabla \varphi dx}  \to 0,\; n\to \infty, \]
and
\[\int\limits_\Omega  {u_t^{{\varepsilon _n}}\varphi dx}  = - \int\limits_\Omega  {\left\langle {{\nu ^{{\varepsilon _n}}},\overrightarrow q  + {\varepsilon _n}I} \right\rangle  \cdot \nabla \varphi dx}. \]
Letting $n\to \infty $, we have
\[\int\limits_\Omega  {{u_t}\varphi dx}  = -\int\limits_\Omega  {\left\langle {\nu ,\overrightarrow q } \right\rangle  \cdot \nabla \varphi dx},\]
i.e.,
\begin{align*}
&u_t=\rm{div} \langle \nu,\overrightarrow{q} \rangle \quad\text{in}\; \mathcal{D}'(\Omega).
\end{align*}
Since $u^{{\varepsilon _n}}\in H^1_0(\Omega)$, then
\[D{{\tilde u}^{{\varepsilon _n}}} = \left\langle {{\nu ^{{\varepsilon _n}}},I} \right\rangle {\mathcal{L}^N}\llcorner\Omega.\]
As
\[D{\tilde u}^{{\varepsilon _n}}\stackrel{*}{\rightharpoondown} D{\tilde u}\quad \text{in}\; \mathcal{M}(\tilde{{\Omega}}),\]
and
\begin{align}
\left\langle {{\nu ^{{\varepsilon _n}}}, I} \right\rangle {\mathcal{L}^N}\llcorner\Omega \stackrel{*}{\rightharpoondown} \left\langle {\nu, I} \right\rangle {\mathcal{L}^N}\llcorner\Omega +\left\langle {{\nu^\infty },{ I }} \right\rangle {\lambda _\nu }\;\text{in}\;\mathcal{M}(\overline{\Omega}),
\end{align}
we obtain
\begin{align}
D{\tilde u}= \left\langle {\nu, I} \right\rangle {\mathcal{L}^N}\llcorner\Omega +\left\langle {{\nu^\infty },{ I }} \right\rangle {\lambda _\nu }\quad\text{in}\;\mathcal{M}(\tilde{{\Omega}}). \label{12}
\end{align}
Therefore,
\begin{align}
Du = \left\langle {\nu,I} \right\rangle {\mathcal{L}^N}\llcorner\Omega + \left\langle {{\nu^\infty },I} \right\rangle {\lambda _\nu }\llcorner \Omega \quad\text{in}\;\mathcal{M}(\Omega),
\end{align}
and by boundary trace theorem (see \cite[Theorem 3.87]{MR1857292}),
\begin{align}
D{\tilde u}\llcorner \partial\Omega =\left\langle {{\nu^\infty },I} \right\rangle {\lambda _\nu }\llcorner \partial\Omega =\left( { - u^\Omega \otimes \vec{n}} \right){\mathcal{H}^{N-1}}\llcorner \partial\Omega  \quad\text{in}\;\mathcal{M}(\overline{\Omega}).
\end{align}

{\textbf{Step 2}}. Prove $$\text{supp }\nu  \subseteq \left\{\Phi = \Phi^{**} \right\}.$$
Set \[{B_\varepsilon }:= \left\{ {A \in {\mathbb{R}^N}\mid {{{\left( {\Phi  + \varepsilon{{\vert I \vert}^2}} \right)}^{**}}(A) = \left( {\Phi  + \varepsilon{{\vert I \vert}^2}} \right)(A)} } \right\},\]
and
\[B := \left\{ {A \in {R^N}: {{\Phi ^{**}}(A) = \Phi (A)} } \right\}.\]
By the definition of $\Phi^{**}$ and ${{\left( {\Phi  + \varepsilon {{\vert I \vert}^2}} \right)}^{**}}$, we arrive at
\[{\Phi ^{**}} + \varepsilon {\vert I \vert^2} \leqslant {\left( {\Phi  + \varepsilon {{\vert I \vert}^2}} \right)^{**}} \leqslant \Phi  + \varepsilon {\vert I \vert^2}.\]
Thus, we get $B\subseteq B_\varepsilon$, for all $\varepsilon>0$.

Similarly, it follows that $B_{\varepsilon_1}\subseteq B_{\varepsilon_2}$, for all $0\leq \varepsilon_1\leq \varepsilon_2$.

Now, assume $\mathop  \cap \limits_{\varepsilon  > 0} {B_\varepsilon } \not\subset B$.
Since $B\subset \mathop  \cap \limits_{\varepsilon  > 0} {B_\varepsilon }$, there exists $x_0 \in \mathop  \cap \limits_{\varepsilon  > 0} {B_\varepsilon }\setminus B$ and $\delta >0$ such that
\begin{align}
\Phi ({x_0}) > {\Phi ^{**}}({x_0}) + 2\delta . \label{se1}
\end{align}
According to Lemma \ref{envelope}, there exist $\alpha_i$, $x_i$, $N$ such that
\[\sum\limits_{i = 1}^{N + 1} {{\alpha _i}{x_i}}  = {x_0},\;\sum\limits_{i = 1}^{N + 1} {{\alpha _i}}  = 1,\;{\alpha _i} \geqslant 0,\]
and
\begin{align}
{\Phi ^{**}}({x_0}) \geqslant \sum\limits_{i = 1}^{N + 1} {{\alpha _i}\Phi ({x_i})}  - \delta . \label{se2}
\end{align}
Observe that
\[ \mathop  \cap \limits_{\varepsilon  > 0} {B_\varepsilon }\subseteq {B_\varepsilon },\;\forall \varepsilon>0,\]
one has
\[{\left( {\Phi  + \varepsilon {{\vert I \vert}^2}} \right)^{**}}({x_0}) = \Phi ({x_0}) + \varepsilon {\vert {{x_0}} \vert^2}.\]
Hence, by \eqref{se1} and \eqref{se2}, it follows that
\[\begin{gathered}
  \sum\limits_{i = 1}^{N + 1} {\left( {{\alpha _i}\Phi ({x_i}) + {\alpha _i}\varepsilon {{\vert {{x_i}} \vert}^2}} \right)}  \hfill \\
   \geqslant {\left( {\Phi  + \varepsilon {{\vert I \vert}^2}} \right)^{**}}({x_0}) \hfill \\
   = \Phi ({x_0}) + \varepsilon {\vert {{x_0}} \vert^2} \hfill \\
   > {\Phi ^{**}}({x_0}) + \varepsilon {\vert {{x_0}} \vert^2} + 2\delta  \hfill \\
   \geqslant \sum\limits_{i = 1}^{N + 1} {{\alpha _i}\Phi ({x_i})}  + \varepsilon {\vert {{x_0}} \vert^2} + \delta. \hfill \\
\end{gathered} \]
Then, we obtain
\[\varepsilon \left( {\sum\limits_{i = 1}^{N + 1} {{\alpha _i}{{\vert {{x_i}} \vert}^2}}  - {{\vert {{x_0}} \vert}^2}} \right) > \delta, \;\forall \varepsilon>0 ,\]
we deduce that it is a contradiction by letting $\varepsilon\rightarrow0$.
\\Thus,
\[\mathop  \cap \limits_{\varepsilon  > 0} {B_\varepsilon } = B.\]
Moreover, recall that \[{\lambda ^{{\varepsilon _n}}} \stackrel{*}{\rightharpoondown}  \lambda \quad\text{in}\;Y\left( {\Omega ;{\mathbb{R}^N}} \right)\]
and
\begin{align*}
\text{supp }\nu^\varepsilon  \subseteq B_{\frac{\varepsilon }{2}}.
\end{align*}
Finally, we get
$$\text{supp }\nu  \subseteq \left\{\Phi = \Phi^{**} \right\}.$$

{\textbf{Step 3.}}
Since
\begin{align*}
&u_t=\rm{div} \langle \nu,\overrightarrow{q} \rangle \quad\text{in}\;\mathcal{D}'(\Omega),
\end{align*}
with ${u_t}\in {L^2}(\Omega )$, we get
\[\left\langle {\nu ,\overrightarrow q } \right\rangle  \in {X_2}(\Omega ).\]
For any  $0\leq\varphi \in C^{\infty}_c(\tilde{\Omega})$, one has
\begin{align}
   &- \int\limits_{\tilde \Omega } {u_t^{{\varepsilon _n}}\varphi {u^{{\varepsilon _n}}}dx}  =  - \int\limits_{\tilde \Omega } {{\rm div}\left( {\left\langle {{\nu ^{{\varepsilon _n}}},\overrightarrow q  + {\varepsilon _n}I} \right\rangle } \right)\varphi {u^{{\varepsilon _n}}}dx} \nonumber \\
   = &\int\limits_{\tilde \Omega } {\varphi \left\langle {{\nu ^{{\varepsilon _n}}},\overrightarrow q  + {\varepsilon _n}I} \right\rangle  \cdot \nabla {u^{{\varepsilon _n}}}dx}  + \int\limits_{\tilde \Omega } {{u^{{\varepsilon _n}}}\left\langle {{\nu ^{{\varepsilon _n}}},\overrightarrow q  + {\varepsilon _n}I} \right\rangle  \cdot \nabla \varphi dx}. \label{st1}
\end{align}
Observe that by \eqref{au2} and \eqref{idp},
\begin{align}
  &\int\limits_{\tilde \Omega } {\varphi \left\langle {{\nu ^{{\varepsilon _n}}},\overrightarrow q  + {\varepsilon _n}I} \right\rangle  \cdot \nabla {u^{{\varepsilon _n}}}dx}  \nonumber\\
   = &\int\limits_{\tilde \Omega } {\varphi \left\langle {{\nu ^{{\varepsilon _n}}},\overrightarrow q  + {\varepsilon _n}I} \right\rangle  \cdot \left\langle {{\nu ^{{\varepsilon _n}}},I} \right\rangle dx}  \nonumber \\
   =& \int\limits_{\tilde \Omega } {\varphi \left\langle {{\nu ^{{\varepsilon _n}}},\overrightarrow q  \cdot I} \right\rangle dx}  + {\varepsilon _n}\int\limits_{\tilde \Omega } {\varphi {{\vert {\left\langle {{\nu ^{{\varepsilon _n}}},I} \right\rangle } \vert}^2}dx}  \nonumber \\
   \geqslant &\int\limits_{\tilde \Omega } {\varphi \left\langle {{\nu ^{{\varepsilon _n}}},\overrightarrow q  \cdot I} \right\rangle dx}.   \label{st2}
\end{align}
From \eqref{cve1},\eqref{cve2},\eqref{cve3} and \eqref{ve1}, we get
\begin{align}
  &\int\limits_{\tilde \Omega } {{u^{{\varepsilon _n}}}\left\langle {{\nu ^{{\varepsilon _n}}},\overrightarrow q  + {\varepsilon _n}I} \right\rangle  \cdot \nabla \varphi dx}  \nonumber\\
   =& \int\limits_{\tilde \Omega } {{u^{{\varepsilon _n}}}\left\langle {{\nu ^{{\varepsilon _n}}},\overrightarrow q } \right\rangle  \cdot \nabla \varphi dx}  + \int\limits_{\tilde \Omega } {{u^{{\varepsilon _n}}}\left( {{\varepsilon _n}\nabla {u^{{\varepsilon _n}}} \cdot \nabla \varphi } \right)dx}  \nonumber \\
   \to& \int\limits_{\tilde \Omega } {u\left\langle {\nu ,\overrightarrow q } \right\rangle  \cdot \nabla \varphi dx} ,\;\;n \to \infty. \label{st3}
\end{align}
On the other hand, by Green's formula of Lemma \ref{green}, it follows that
\begin{align}
  & - \int\limits_{\tilde \Omega } {{u_t}\varphi udx}  \nonumber\\
   =&  - \int\limits_{\tilde \Omega } {{\rm div}\left( {\left\langle {\nu ,\overrightarrow q } \right\rangle } \right)\varphi udx}  \nonumber \\
   =& \int\limits_{\tilde \Omega } {\varphi \left( {\left\langle {\nu ,\overrightarrow q } \right\rangle ,D\tilde u} \right)}  + \int\limits_{\tilde \Omega } {u\left\langle {\nu ,\overrightarrow q } \right\rangle  \cdot \nabla \varphi dx}  \nonumber \\
   =& \int\limits_{\tilde \Omega } \varphi  \left\langle {\nu ,\overrightarrow q } \right\rangle  \cdot \left( {\left\langle {\nu,I} \right\rangle dx + \left\langle {{\nu^\infty },I} \right\rangle d{\lambda _\nu }} \right) + \int\limits_{\tilde \Omega } {u\left\langle {\nu ,\overrightarrow q } \right\rangle  \cdot \nabla \varphi dx}.    \label{st4}
\end{align}
Hence, by \eqref{st1},\eqref{st2},\eqref{st3} and \eqref{st4}, letting $n\to \infty$, we have
\begin{align*}
  &\int\limits_{\tilde \Omega } \varphi  \left\langle {\nu ,\overrightarrow q } \right\rangle  \cdot \left( {\left\langle {\nu,I} \right\rangle dx + \left\langle {{\nu^\infty },I} \right\rangle d{\lambda _\nu }} \right)  \\
   \geqslant &\int\limits_{\tilde \Omega } \varphi  \left( {\left\langle {\nu,\overrightarrow q  \cdot I} \right\rangle dx + \left\langle {{\nu^\infty },{{(\overrightarrow q  \cdot I)}^\infty }} \right\rangle d{\lambda _\nu }} \right).
\end{align*}
Furthermore, by the arbitrariness of $\phi$ and the definition of $\tilde{u}$,
\begin{align}
&\left\langle {\nu ,\overrightarrow q } \right\rangle  \cdot \left(  \left\langle {\nu,I} \right\rangle {\mathcal{L}^N}\llcorner\Omega + \left\langle {{\nu^\infty },I} \right\rangle {\lambda _\nu }\llcorner \overline{\Omega}\right)  \nonumber
\\ \geq &\left\langle {\nu,\overrightarrow{q}\cdot I} \right\rangle {\mathcal{L}^N}\llcorner\Omega + \left\langle {{\nu^\infty },(\overrightarrow{q}\cdot I)^\infty} \right\rangle {\lambda _\nu }\llcorner \overline{\Omega}        \label{15}
\end{align}
holds in $\mathcal{M}(\overline{\Omega})$. $\hfill\blacksquare$

\begin{remark}
Actually, the {\bf{(JF)}} inequality \eqref{jf} is also an equality. According to monotonicity of $\overrightarrow{p}$, one has
$(\overrightarrow{p}(A)-\overrightarrow{p}(A_0))\cdot(A-A_0)\geq 0$. By \eqref{supp}, taking $A_0=\langle \nu, I\rangle$, we get
$\langle \nu, \overrightarrow{q}\cdot I\rangle \geq \langle \nu, \overrightarrow{q}\rangle\cdot \langle \nu, I\rangle.$
On the other hand, by $\rm (SH_4)$, we have
$\langle \nu ^\infty , (\overrightarrow{q}\cdot I)^\infty \rangle \geq \langle \nu , \overrightarrow{q}\rangle\cdot \langle \nu^\infty, I\rangle.$
These two inequalities lead to the conclusion.
\end{remark}

\section{Uniqueness Results}
The {\bf{(JF)}} inequality \eqref{jf} describes a property of the solution upon which the uniqueness proof relies. Similar to the situation of Young measure solutions, we note that there is no claim that the parametrized measure $\lambda$ is unique, which is false in general.
\begin{theorem} \label{unique}
Under the assumptions of {\bf{(SH)}}. If both $(u_1,\lambda_1)$ and $(u_2,\lambda_2)$ are the generalized Young measure solutions of $(\mathcal{P})$, then
\[u_1=u_2\quad\textrm{{a.e. in }} \Omega_T\]
for every $T>0$.
\end{theorem}
{\bf{Proof}}: Note $\lambda_1=(\nu_1,\lambda^1_\nu,\nu^\infty _1)$, $\lambda_2=(\nu_2,\lambda^2_\nu,\nu^\infty _2)$.
Since
\begin{align*}
    &-\int_0^s {\int\limits_{\tilde \Omega } {\frac{{\partial ({u_1} - {u_2})}}{{\partial t}}({u_1} - {u_2})dx} } dt \\
   = &- \int_0^s {\int\limits_{\tilde \Omega } {{\rm div}\left( {\left\langle {{\nu _1},\overrightarrow q } \right\rangle  - \left\langle {{\nu _2},\overrightarrow q } \right\rangle } \right)({u_1} - {u_2})dx} } dt\\
   = &\int_0^s {\int\limits_{\tilde \Omega } {\left( {\left\langle {{\nu _1},\overrightarrow q } \right\rangle ,D{{\tilde u}_1}} \right)} } dt - \int_0^s {\int\limits_{\tilde \Omega } {\left( {\left\langle {{\nu _2},\overrightarrow q } \right\rangle ,D{{\tilde u}_1}} \right)} } dt   \\&+
  \int_0^s {\int\limits_{\tilde \Omega } {\left( {\left\langle {{\nu _2},\overrightarrow q } \right\rangle ,D{{\tilde u}_2}} \right)} } dt - \int_0^s {\int\limits_{\tilde \Omega } {\left( {\left\langle {{\nu _1},\overrightarrow q } \right\rangle ,D{{\tilde u}_2}} \right)} } dt\\
   = & \int_0^s {\int\limits_{\tilde \Omega } {\left\langle {{\nu _1},\overrightarrow q } \right\rangle \left( {\left\langle {{\nu_1},I} \right\rangle dx + \left\langle {\nu_1^\infty ,I} \right\rangle d\lambda _\nu ^1} \right)} } dt  \\&+
  \int_0^s {\int\limits_{\tilde \Omega } {\left\langle {{\nu _2},\overrightarrow q } \right\rangle \left( {\left\langle {{\nu_2},I} \right\rangle dx + \left\langle {\nu_2^\infty ,I} \right\rangle d\lambda _\nu ^2} \right)} } dt  \\&-
  \int_0^s {\int\limits_{\tilde \Omega } {\left\langle {{\nu _2},\overrightarrow q } \right\rangle \left( {\left\langle {{\nu_1},I} \right\rangle dx + \left\langle {\nu_1^\infty ,I} \right\rangle d\lambda _\nu ^1} \right)} } dt
  \\&-\int_0^s {\int\limits_{\tilde \Omega } {\left\langle {{\nu _1},\overrightarrow q } \right\rangle \left( {\left\langle {{\nu_2},I} \right\rangle dx + \left\langle {\nu_2^\infty ,I} \right\rangle d\lambda _\nu ^2} \right)} } dt \\
   \geqslant & \int_0^s {\int\limits_{\tilde \Omega } {\left( {\left\langle {{\nu _1},\overrightarrow q  \cdot I} \right\rangle dx + \left\langle {\nu_1^\infty ,{{(\overrightarrow q  \cdot I)}^\infty }} \right\rangle d\lambda _\nu ^1} \right)} } dt  \\&+
  \int_0^s {\int\limits_{\tilde \Omega } {\left( {\left\langle {{\nu _2},\overrightarrow q  \cdot I} \right\rangle dx + \left\langle {\nu_2^\infty ,{{(\overrightarrow q  \cdot I)}^\infty }} \right\rangle d\lambda _\nu ^2} \right)} } dt \\&-
  \int_0^s {\int\limits_{\tilde \Omega } {\left\langle {{\nu _2},\overrightarrow q } \right\rangle \left( {\left\langle {{\nu_1},I} \right\rangle dx + \left\langle {\nu_1^\infty ,I} \right\rangle d\lambda _\nu ^1} \right)} } dt  \\&-
  \int_0^s {\int\limits_{\tilde \Omega } {\left\langle {{\nu _1},\overrightarrow q } \right\rangle \left( {\left\langle {{\nu_2},I} \right\rangle dx + \left\langle {\nu_2^\infty ,I} \right\rangle d\lambda _\nu ^2} \right)} } dt\\
   = & \int_0^s {\int\limits_{\tilde \Omega } {\left( {\left\langle {{\nu _1},\overrightarrow q  \cdot I} \right\rangle  + \left\langle {{\nu _2},\overrightarrow q  \cdot I} \right\rangle  - \left\langle {{\nu _2},\overrightarrow q } \right\rangle \left\langle {{\nu_1},I} \right\rangle  - \left\langle {{\nu _1},\overrightarrow q } \right\rangle \left\langle {{\nu_2},I} \right\rangle } \right)dx} } dt  \\&+
  \int_0^s {\int\limits_{\tilde \Omega } {\left( {\left\langle {\nu_1^\infty ,{{(\overrightarrow q  \cdot I)}^\infty }} \right\rangle  - \left\langle {{\nu _2},\overrightarrow q } \right\rangle \left\langle {\nu_1^\infty ,I} \right\rangle } \right)d\lambda _\nu ^1} } dt  \\&+
   \int_0^s {\int\limits_{\tilde \Omega } {\left( {\left\langle {\nu_2^\infty ,{{(\overrightarrow q  \cdot I)}^\infty }} \right\rangle  - \left\langle {{\nu _1},\overrightarrow q } \right\rangle \left\langle {\nu_2^\infty ,I} \right\rangle } \right)d\lambda _\nu ^2} } dt  \\
  : = & {I_1} + {I_2} + {I_3}.
\end{align*}
Observe that by \eqref{supp},
\begin{align*}
  &\left\langle {{\nu _1},\overrightarrow q  \cdot I} \right\rangle  + \left\langle {{\nu _2},\overrightarrow q  \cdot I} \right\rangle  - \left\langle {{\nu _2},\overrightarrow q } \right\rangle \left\langle {{\nu_1},I} \right\rangle  - \left\langle {{\nu _1},\overrightarrow q } \right\rangle \left\langle {{\nu_2},I} \right\rangle \\
   = &\int\limits_{{\mathbb{R}^N}} {\int\limits_{{\mathbb{R}^N}} {\overrightarrow q (A) \cdot I(A) + \overrightarrow q (B) \cdot I(B) - \overrightarrow q (A) \cdot I(B) - \overrightarrow q (B) \cdot I(A)d{\nu _1}(A)d{\nu _2}(B)} }  \\
   = &\int_{\left\{ {\Phi  = {\Phi ^{**}}} \right\}} {\int_{\left\{ {\Phi  = {\Phi ^{**}}} \right\}} {\left( {\overrightarrow q (A) - \overrightarrow q (B)} \right) \cdot \left( {I(A) - I(B)} \right)d{\nu _1}(A)d{\nu _2}(B)} }  \\
   = &\int_{\left\{ {\Phi  = {\Phi ^{**}}} \right\}} {\int_{\left\{ {\Phi  = {\Phi ^{**}}} \right\}} {\left( {\overrightarrow p (A) - \overrightarrow p (B)} \right) \cdot \left( {I(A) - I(B)} \right)d{\nu _1}(A)d{\nu _2}(B)} }  \\
   \geqslant & 0.
\end{align*}
Hence, it follows that $I_1\geq 0$.

On the other hand, by $\rm (SH_4)$ and \eqref{supp},
\begin{align*}
  &\left\langle {\nu_1^\infty ,{{(\overrightarrow q  \cdot I)}^\infty }} \right\rangle  - \left\langle {{\nu _2},\overrightarrow q } \right\rangle \left\langle {\nu_1^\infty ,I} \right\rangle  \\
   = &\int_{\left\{ {\Phi  = {\Phi ^{**}}} \right\}} {\int_{\partial {\mathbb{B}^N}} {{{(\overrightarrow q  \cdot I)}^\infty }(A) - \overrightarrow q (B) \cdot I(A)d\nu_1^\infty (A)d{\nu _2}(B)} }   \\
   \geqslant & 0.
\end{align*}
Therefore, $I_2\geq 0$. Similarly, we get $I_3\geq 0$.

Now, by
\[\int_0^s {\int\limits_{\tilde \Omega } {\frac{{\partial ({u_1} - {u_2})}}{{\partial t}}({u_1} - {u_2})dx} } dt\leq 0,\]
we have
\[u_1=u_2\quad \text{a.e. in}\; \Omega, \; \forall T>0,\]
and this leads to the uniqueness. $\hfill\blacksquare$

\section{Equivalence between Generalized Young Measure Solutions and Strong Solutions}

It doesn't make much sense just to propose a concept of generalized solution with existence and uniqueness. Next, we will establish the relationship between this generalized solution and other types of generalized solutions.

When the parabolic variational integral $\Phi$ is convex, Andreu et al.\cite{2004Parabolic} combined nonlinear semigroup theory with Anzellotti's dual theory\cite{Gabriele1983Pairings} to give a definition of semigroup strong solutions to this equation, and systematically studied the existence, uniqueness and asymptotic behavior of this kind of solutions.
In this section, we study the relationship between generalized Young measure solutions and strong solutions (\cite[Definition 6.5]{2004Parabolic}).

The additional assumptions {\bf{(H)}} are the following:
\\$\rm(H_1)$ $\widehat{\Phi} (\xi,s)$ is continuous on $\mathbb{R}^N \times [0,+\infty)$ and convex in $(\xi,s)$, where
\[\widehat \Phi (\xi ,s): = \left\{ \begin{gathered}
  \Phi \left( {\frac{\xi }{s}} \right)s,\quad if\;s > 0, \hfill \\
  {\Phi ^\infty }\left( \xi  \right),\quad if\;s = 0. \hfill \\
\end{gathered}  \right.\]
$\rm(H_2)$ ${\Phi ^\infty }\left( { - \xi } \right) = {\Phi ^\infty }\left( \xi  \right)$  and  $\overrightarrow q (\xi ) \cdot \xi  \geqslant 0$ holds, $\forall \xi \in \mathbb{R}^N$.

\begin{definition}\label{sms}
A measurable function $u:(0,T)\times \Omega \to \mathbb{R}$ is a strong solution of $(\mathcal{P})$ in $\Omega_T$ if $u\in C([0,T],L^2(\Omega))$, $u(0)=u_0$, $u'(t)\in L^2(\Omega)$, $u(t)\in BV(\Omega)\cap L^2(\Omega)$, $\overrightarrow{q}(\nabla u(t))\in X(\Omega)_1$ a.e.\;$t\in[0,T]$, and for almost all $t\in[0,T]$ $u(t)$ satisfies:
\begin{align}
&u'(t)={\rm{div}}(\overrightarrow{q}(\nabla u(t))) \;\;\text{in}\; \mathcal{D}(\Omega), \label{d51}
\\&\overrightarrow{q}(\nabla u(t)) \cdot D^s u(t)=\Phi^\infty (D^s u(t)),\label{d52}
\\&[\overrightarrow{q}(\nabla u(t)), \vec{n}]\in {\rm sign}(-u^\Omega (t))\Phi^\infty (\vec{n})\quad \mathcal{H}^{N-1}\text{-a.e.}\; \text{on}\; \partial \Omega.\label{d53}
\end{align}
\end{definition}

Using the theory of nonlinear semigroups, Andreu et al. proved the existence and uniqueness of strong solutions \cite[Theorem 6.6]{2004Parabolic} and the following lemma follows.
\begin{lemma}
Under the assumptions of $\bf{(SH)}$ and $\bf{(H)}$. Given $u_0 \in L^2(\Omega)$, there exists a unique strong solution $u$ of $(\mathcal{P})$ in $\Omega_T$ for every $T>0$ such that $u(0)=u_0$.
\end{lemma}

Assume that $\Phi$ satisfies $\bf(SH)$ and $u$ is a solution in the Definition \ref{sms}, by calculation, we can verify that $(u,\sigma_{Du})$ is a generalized gradient Young solution. The following theorem follows.
\begin{theorem}
If $u$ is a strong solution of $(\mathcal{P})$ in the Definition \ref{sms}, then $(u,\sigma_{Du})$ is a generalized gradient Young measure solution of $(\mathcal{P})$.
\end{theorem}

On the other hand, if $\Phi$ is convex and satisfies {\bf{(SH)}}, the generalized Young solution of $(\mathcal{P})$ uniquely determines a strong solution in Definition \ref{sms}.

\begin{theorem}
Let $u_0\in BV(\Omega)\cap L^{\infty}(\Omega)$. If $(u,\lambda)$ is a generalized Young measure solution of $(\mathcal{P})$, then $u$ is a strong solution of $(\mathcal{P})$ in Definition \ref{sms}.
\end{theorem}
{\bf{Proof}}:
Since $\Phi$ is convex function, by the usual theory of monotone operators, there exists a unique weak solution (of usual sense) to problem $(\mathcal{P}_\varepsilon)$,
then $(u^{\varepsilon}, \delta_{\nabla u_{\varepsilon}})$ is a Young measure solution to the problem $(\mathcal{P}_\varepsilon)$.

Using $(u^{\varepsilon}, \delta_{\nabla u_{\varepsilon}})$, by the method in Theorem \ref{exist}, there exist a subsequence $(u^{\varepsilon_n}, \delta_{\nabla u_{\varepsilon_n}})$ and a generalized gradient Young measure solution $(\bar{u},\bar{\lambda})$ ($\bar{\lambda}:=(\bar{\nu},\bar{\lambda}_{\nu},\bar{\nu}^\infty)$) satisfying that
\[\begin{gathered}
  {u^{{\varepsilon _n}}} \to \bar u\quad{\text{in}}\;{L^2}(\Omega ), \hfill \\
  {\varepsilon _n}\vert {\nabla {u^{{\varepsilon _n}}}} \vert \to 0\quad{\text{in}}\;{L^2}(\Omega ), \hfill \\
  \left\langle {{\delta _{\nabla {u^{{\varepsilon _n}}}}},\overrightarrow{q} } \right\rangle \mathop  \rightharpoonup \limits^* \left\langle {\bar \nu ,\overrightarrow{q}} \right\rangle \quad{\text{in}}\;{L^\infty }(\Omega ;{\mathbb{R}^N}), \hfill \\
u_t^{{\varepsilon _n}} {\rightharpoonup} {\bar{u}_t}\quad\text{in}\;{L^2}(\Omega ). \hfill \\
\end{gathered} \]
Since
\begin{align}
{{\bar u}_t} = {\text{div}}\langle \bar \nu ,\overrightarrow{q}\rangle \quad{\text{in}}\;\mathcal{D}'(\Omega ), \label{d56}
\end{align}
we get $\langle \bar \nu ,\overrightarrow{q}\rangle \in X_2(\Omega) $.

Let $0\leq\theta \in \mathcal{D}(\Omega)$ and $\eta\in C^1(\overline{\Omega})$. By the convexity of $\Phi$ and $\overrightarrow{q}=\nabla \Phi$,
it follows that
\begin{align}
\int\limits_\Omega  {\theta \left( {\overrightarrow{q}_{{\varepsilon _n}}\left( {\nabla {u^{{\varepsilon _n}}}} \right)- \overrightarrow{q}_{\varepsilon _n}\left( {\nabla \eta } \right)}\right) \cdot \left( {\nabla {u^{{\varepsilon _n}}} - \nabla \eta } \right)dx}  \geqslant 0. \label{mono}
\end{align}
Since
\begin{align*}
 &\int\limits_\Omega  {\theta \overrightarrow{q}_{\varepsilon _n} \left( {\nabla {u^{{\varepsilon _n}}}} \right) \cdot \nabla \left( {{u^{{\varepsilon _n}}} - \eta } \right)dx}
  \\= &- \int\limits_\Omega  {\theta u_t^{{\varepsilon _n}} \left( {{u^{{\varepsilon _n}}} - \eta } \right)dx}
  -\int\limits_\Omega  {\left( {{u^{{\varepsilon _n}}} - \eta } \right)\nabla \theta  \cdot \overrightarrow{q}_{\varepsilon _n}\left( {\nabla {u^{{\varepsilon _n}}}} \right)dx},
\end{align*}
we have
\begin{align*}
  &\mathop {\lim }\limits_{n \to \infty } \int\limits_\Omega  {\theta \overrightarrow{q}_{\varepsilon _n}\left( {\nabla {u^{{\varepsilon _n}}}} \right) \cdot \nabla \left( {{u^{{\varepsilon _n}}} - \eta } \right)dx}   \\
  =  &  - \int\limits_\Omega  {\theta \bar{u} _t \left( {\bar{u} - \eta } \right)dx} -\int\limits_\Omega  {\left( {\bar{u} - \eta } \right)\nabla \theta  \cdot \left\langle {\bar \nu ,\overrightarrow{q}} \right\rangle dx}                             \\
   = & - \int\limits_\Omega  {\theta {\rm{div}}\left( {\left\langle {\bar \nu ,\overrightarrow{q}} \right\rangle } \right)\left( {\bar{u} - \eta } \right)dx}  - \int\limits_\Omega  {\left( {\bar{u} - \eta } \right)\nabla \theta  \cdot \left\langle {\bar \nu ,\overrightarrow{q}} \right\rangle dx}  \\
   =& \int\limits_\Omega  {\theta \left( {\left\langle {\bar \nu ,\overrightarrow{q}} \right\rangle, D\left( {\bar u - \eta } \right)} \right)}.
\end{align*}
From this, and
\[\mathop {\lim }\limits_{n \to \infty } \int\limits_\Omega  {\theta \overrightarrow{q}_{\varepsilon _n}\left( {\nabla \eta } \right) \cdot \nabla \left( {{u^{{\varepsilon _n}}} - \eta } \right)dx}  = \int\limits_\Omega  {\theta \overrightarrow{q}\left( {\nabla \eta } \right)\cdot D\left( {\bar u - \eta } \right)} ,\]
letting $n\to \infty$ in \eqref{mono}, we obtain
\[\int\limits_\Omega  {\theta \left( {\left\langle {\bar \nu ,\overrightarrow{q}} \right\rangle  - \overrightarrow{q}\left( {\nabla \eta } \right),D\left( {\bar u - \eta } \right)} \right)}  \geqslant 0.\]
Hence, \[{\left( {\left\langle {\bar \nu ,\overrightarrow{q}} \right\rangle  - \overrightarrow{q}\left( {\nabla \eta } \right),D\left( {\bar u - \eta } \right)} \right)}\geqslant 0 \quad \text{in } \mathcal{M}(\Omega),\]
then, its absolutely continuous part satisfies that
\[\left( {\left\langle {\bar \nu ,\overrightarrow{q}} \right\rangle  - \overrightarrow{q}\left( {\nabla \eta } \right)} \right) \cdot \nabla \left( {\bar u - \eta } \right) \geqslant 0 \quad \text{a.e. in }\Omega .\]

Since $C^1(\overline{\Omega})$ is separable, by taking a countable set dense in $C^1(\overline{\Omega})$, we have that the above inequality holds for all $x\in \Omega'$, where $\Omega'\subseteq \Omega$ is such that $\mathcal{L}^N(\Omega \setminus \Omega')=0$, and all $\eta\in C^1(\overline{\Omega})$.
Now, fixed $x\in \Omega'$ and given $A\in \mathbb{R}^N$, there is $\eta\in C^1(\overline{\Omega})$ such that $\nabla \eta =A$.
Thus,
\[\left( {\left\langle {\bar \nu ,\overrightarrow{q}} \right\rangle - \overrightarrow{q}\left( A \right)} \right) \cdot \left( {\nabla \bar u(x) - A} \right) \geqslant 0,\;\forall A \in \mathbb{R}^N.\]
By choosing $A=\nabla \bar{u}(x)\pm \epsilon \xi$, $\forall \xi \in \mathbb{R}^N$, and letting $\epsilon \to 0^+$, we obtain
\begin{align}
\left\langle {\bar \nu ,\overrightarrow{q}} \right\rangle  = \overrightarrow{q}(\nabla {\bar u}(x))\quad\text{a.e.} \;x \in \Omega . \label{eql}
\end{align}

Similar to \eqref{12} and \eqref{15}, we have
\begin{align}
D{\tilde{\bar{u}}}= \left\langle {\bar{\nu}, I} \right\rangle {\mathcal{L}^N}\llcorner\Omega +\left\langle {{\bar{\nu}^\infty },{I }} \right\rangle {\bar{\lambda} _\nu }\quad\text{in}\;\mathcal{M}(\tilde{{\Omega}}), \label{54}
\end{align}
and
\begin{align}
&\left\langle {\bar{\nu} ,\overrightarrow q } \right\rangle  \cdot \left(  \left\langle {\bar{\nu},I} \right\rangle {\mathcal{L}^N}\llcorner\Omega + \left\langle {{\bar{\nu}^\infty },I} \right\rangle {\bar{\lambda} _\nu }\llcorner \overline{\Omega}\right) \nonumber\\ \geqslant &
\left(  \left\langle {\bar{\nu},\overrightarrow{q}\cdot I} \right\rangle {\mathcal{L}^N}\llcorner\Omega + \left\langle {{\bar{\nu}^\infty },(\overrightarrow{q}\cdot I)^\infty} \right\rangle {\bar{\lambda} _\nu }\llcorner \overline{\Omega}\right),
\end{align}
in $\mathcal{M}(\overline{\Omega})$.

Thus,
\begin{align}
&\left( {\left\langle {\bar \nu ,\overrightarrow{q}} \right\rangle \left\langle {\bar \nu ,I} \right\rangle  - \left\langle {\bar \nu ,\overrightarrow{q} \cdot I} \right\rangle } \right){\mathcal{L}^N} \llcorner \Omega  \nonumber\\&+ \left( {\left\langle {\bar \nu ,\overrightarrow{q}} \right\rangle \left\langle {{{\bar \nu }^\infty },I} \right\rangle  - \left\langle {{{\bar \nu }^\infty },{{(\overrightarrow{q} \cdot I)}^\infty }} \right\rangle } \right){{\bar \lambda }_\nu } \llcorner \overline{\Omega}  \geqslant 0,
\end{align}
and its singular parts satisfies that
\[\left( {\left\langle {\bar \nu ,\overrightarrow{q}} \right\rangle \left\langle {{{\bar \nu }^\infty },I} \right\rangle  - \left\langle {{{\bar \nu }^\infty },{{(\overrightarrow{q} \cdot I)}^\infty }} \right\rangle } \right)\bar \lambda _\nu ^s \llcorner \overline{\Omega}  \geqslant 0\quad\text{in}\;\mathcal{M}(\overline{\Omega}).\]
Then, it follows that
\begin{align}
\left\langle {\bar \nu ,\overrightarrow{q}} \right\rangle \left\langle {{{\bar \nu }^\infty },I} \right\rangle  \geqslant \left\langle {{{\bar \nu }^\infty },{{(\overrightarrow{q} \cdot I)}^\infty }} \right\rangle \quad\bar \lambda _\nu ^s {\text{-a.e.}}\;{\text{in}}\;\overline{\Omega}. \label{ne1}
\end{align}
By \eqref{54}, we get
\begin{align}
{D^s}\tilde {\bar u} = \left\langle {{{\bar \nu }^\infty },{{(\overrightarrow{q} \cdot I)}^\infty }} \right\rangle \bar \lambda _\nu ^s\quad{\text{in}}\;\mathcal{M}(\overline{\Omega} ). \label{57}
\end{align}
Since $\bar \lambda$ is a generalized gradient Young measure, by $(ii)$ of Lemma \ref{charact}, it follows that
\[\left\langle {{{\bar \nu }^\infty },{{(\overrightarrow{q} \cdot I)}^\infty }} \right\rangle  = \left\langle {{{\bar \nu }^\infty },{\Phi ^\infty }} \right\rangle  \geqslant {\Phi ^\infty }\left( {\left\langle {{{\bar \nu }^\infty },I} \right\rangle } \right)\quad\bar \lambda _\nu ^s {\text{-a.e.}}\;{\text{in}}\;\overline{\Omega},\]
i.e.,
\begin{align}
\left\langle {{{\bar \nu }^\infty },{{(\overrightarrow{q} \cdot I)}^\infty }} \right\rangle \bar \lambda _\nu ^s \geqslant {\Phi ^\infty }\left( {\left\langle {{{\bar \nu }^\infty },I} \right\rangle } \right)\bar \lambda _\nu ^s = {\Phi ^\infty }\left( {\frac{{{D^s}\tilde {\bar u}}}{{\vert {{D^s}\tilde {\bar u}} \vert}}} \right)\vert {{D^s}\tilde {\bar u}} \vert\quad{\text{in}}\;\mathcal{M}(\overline{\Omega} ). \label{58}
\end{align}
On the other hand, by \eqref{57}, we have
\[\left\langle {\bar \nu ,\overrightarrow{q}} \right\rangle \left\langle {{{\bar \nu }^\infty },I} \right\rangle \bar \lambda _\nu ^s = \left\langle {\bar \nu ,\overrightarrow{q}} \right\rangle {D^s}\tilde {\bar u }= \left\langle {\bar \nu ,\overrightarrow{q}} \right\rangle \left\langle {{\delta _p},I} \right\rangle \vert {{D^s}\tilde {\bar u}} \vert,\]
where $p = \frac{{{D^s}\tilde {\bar u}}}{{\vert {{D^s}\tilde {\bar u}} \vert}}$.

Then, by $\rm (SH_4)$, we find
\begin{align*}
 & \left\langle {\bar \nu ,\overrightarrow{q}} \right\rangle \left\langle {{\delta _p},I} \right\rangle   \\
   = &\int\limits_{{\mathbb{R}^N}} {\int\limits_{\partial {\mathbb{B}^N}} {\overrightarrow{q}(B) \cdot Ad{\delta _p}(A)d\bar \nu (B)} }  \\
   \leqslant & \int\limits_{{\mathbb{R}^N}} {\int\limits_{\partial {\mathbb{B}^N}} {{{(\overrightarrow{q} \cdot I)}^\infty }(A)d{\delta _p}(A)d\bar \nu (B)} }  \\
   = &{\Phi ^\infty }\left( {\frac{{{D^s}\tilde {\bar u}}}{{\vert {{D^s}\tilde {\bar u}} \vert}}} \right)
\end{align*}
holds $\vert {{D^s}\tilde {\bar u}} \vert {\text{-a.e.}}\;{\text{in}}\;\overline{\Omega}$.

Combining this inequality with \eqref{ne1} and \eqref{58}, we obtain
\[\left\langle {\bar \nu ,\overrightarrow{q}} \right\rangle {D^s}\tilde {\bar u} = {\Phi ^\infty }\left( {\frac{{{D^s}\tilde {\bar u}}}{{\vert {{D^s}\tilde {\bar u}} \vert}}} \right)\vert {{D^s}\tilde {\bar u}} \vert = {\Phi ^\infty }\left( {{D^s}\tilde {\bar u}} \right) \quad{\text{in}}\;\mathcal{M}(\overline{\Omega} ).\]
Therefore,
\begin{align}
\left\langle {\bar \nu ,\overrightarrow{q}} \right\rangle {{D^s}\bar{u}} \;\llcorner \Omega  = {\Phi ^\infty }\left( {{D^s}\bar u} \right) \llcorner \Omega \quad{\text{in}}\;\mathcal{M}(\overline{\Omega} ), \label{d54}
\end{align}
and
\begin{align}
\left\langle {\bar \nu ,\overrightarrow{q}} \right\rangle {D^s}\tilde {\bar u} \llcorner \partial \Omega  = {\Phi ^\infty }\left( {{D^s}\tilde {\bar u}} \right)\quad{\text{in}}\;\mathcal{M}(\overline{\Omega} ).\label{d55}
\end{align}
Since $\bar{u}\in BV(\Omega)$, we get
\begin{align}
{D^s}\tilde {\bar u} \llcorner \partial \Omega  = D\tilde {\bar u} \llcorner \partial \Omega  = \left( { - {{\bar u}^\Omega } \otimes \vec n} \right){\mathcal{H}^{N - 1}} \llcorner \partial \Omega \quad{\text{in}}\;\mathcal{M}(\overline{\Omega} ). \label{d57}
\end{align}
 According to Theorem \ref{unique}, $\bar{u}=u$.

Combining \eqref{d56},\eqref{eql} and \eqref{d54}--\eqref{d57}, we arrive at \eqref{d51}--\eqref{d53}. The proof is complete.  $\hfill\blacksquare$

\bmhead{Acknowledgments}


\bmhead{Data Availability}

The authors declare that data sharing not applicable to this article, as no datasets were generated or analysed during the study.

%

%
%
%
%
%
%
%
%

\bibliography{bibliography}


\end{document}